


\magnification=\magstep1
\baselineskip=14pt  
\abovedisplayskip=10pt plus 2pt minus 2pt
\belowdisplayskip=10pt plus 2pt minus 2pt
\font\bigbf=cmbx10 scaled \magstep1  
\def\prf{\noindent {\bf Proof.\ \ }}  
\def\square{\smash{\raise8pt\hbox{\leavevmode\hbox{\vrule\vtop{\vbox
	{\hrule width9pt}\kern9pt\hrule width9pt}\vrule}}}}
\def\qed{\hfill\square}  


\def\res{{\rm Res}} 
\def\P{{\bf P}}     
\def\C{{\bf C}}     
\def\R{{\bf R}}     
\def\Q{{\bf Q}}     
\def\Z{{\bf Z}}     
\def\ox{{\cal O}_X} 
\def\U {{\cal U}}   
\def\ip{\,\smash{\raise7pt\hbox{\leavevmode\hbox{\vtop{\kern7pt\hrule
  width4pt}\vrule}}}\,\,}  
\def\mapname#1{\ \smash{\mathop{\longrightarrow}\limits^{#1}}\ } 
\def\dbar{\overline{\partial}} 

\centerline{\bigbf Toric Residues}
\medskip
\medskip
\centerline{David A. Cox}
\centerline{Department of Mathematics and Computer Science}
\centerline{Amherst College}
\centerline{Amherst, MA 01002}
\centerline{dac@cs.amherst.edu}
\medskip

	The Grothendieck local residue symbol
$$\res_{0}\left( {g\, dx_0\wedge \cdots\wedge dx_n \over f_0\cdots f_n}
	\right) = {1 \over (2\pi i)^{n+1}} \int_{|f_i| =
	\epsilon} {g\, dx_0\wedge \cdots\wedge dx_n \over f_0\cdots f_n}
	\leqno(1)$$
(see [GH, Chapter 5]) is defined whenever $g,f_0,\dots,f_n$ are
holomorphic in a neighborhood of $0 \in \C^{n+1}$ and $f_0,\dots,f_n$
don't vanish simultaneously except at $0$.  In [PS, 12.10], it was
observed that when $f_0,\dots,f_n$ are homogeneous of degree $d$ and
$g$ is homogeneous of degree $\rho = (n+1)(d-1)$, the residue symbol
has the following nice properties.

\proclaim Quotient Property. The map
$$ g \mapsto \res_{0}\left( {g\, dx_0\wedge \cdots\wedge dx_n \over
	f_0\cdots f_n} \right) $$
induces an isomorphism
$$\C[x_0,\dots,x_n]_\rho/\langle f_0,\dots,f_n\rangle_\rho \simeq \C$$
(the subscript refers to the graded piece in degree $\rho$) uniquely
characterized by the fact that the Jacobian determinant $J =
\det(\partial f_i/\partial x_j)$ maps to $d^{n+1}$.

\proclaim Trace Property. \v Cech cohomology gives a naturally
defined cohomology class $[\omega_g] \in H^n(\P^n,\Omega^n_{\P^n})$
such that under the trace map ${\rm Tr}_{\P^n} : H^n(\P^n,\Omega^n_{\P^n})
\simeq \C$, we have
$$\res_{0}\left( {g\, dx_0\wedge \cdots\wedge dx_n \over f_0\cdots f_n}
	\right) = {\rm Tr}_{\P^n}([\omega_g])\ . $$
(We define $[\omega_g]$ in \S1.)

	In this paper, we will show how these properties of residues
can be generalized to an arbitrary projective toric variety.  The
paper is organized into six sections as follows.  In \S1, we define
the cohomology class $[\omega_g] \in H^n(\P^n,\Omega^n_{\P^n})$, and
then \S2 generalizes this to define toric residues in terms of a toric
analog of the Trace Property.  We recall some commutative algebra
associated with toric varieties in \S3, and \S4 introduces a toric
version of the Jacobian, which is defined as the determinant of a
certain exact sequence.  In \S5, we show that the toric residue is
uniquely characterized using a toric analog of the Quotient Property.
Then \S6 explores different ways of representing the toric residue as
an integral, and an appendix discusses the relation between the trace
map and the Dolbeault isomorphism.

	This paper had its origin in a series of conversations with
David Morrison, and we are grateful to Morrison for asking the
question that led to the results presented here.  We also would like
to thank Eduardo Cattani and Alicia Dickenstein for their comments on
an earlier version of the manuscript.  The research for this paper
was supported by NSF grant DMS-9301161.
\medskip

	{\bf \S1. The Definition of $\omega_g$ for $\P^n$.}  Suppose
that $f_0,\dots,f_n$ are homogeneous polynomials of degree $d$
which don't vanish simultaneously on $\C^{n+1}$ except at the origin,
and let $g$ be homogeneous of degree $\rho = (n+1)(d-1)$.  Then
consider the $n$-form
$$\Omega = \sum_{i=0}^n (-1)^i x_i\, dx_0\wedge \cdots \wedge
	\widehat{dx_i}\wedge \cdots \wedge dx_n\ .\leqno(2)$$
As is well-known (see [G, \S2]), our assumptions on $g$ and
$f_0,\dots,f_n$ imply that
$$\omega_g = {g\,\Omega \over f_0\cdots f_n}$$
descends to a meromorphic $n$-form on $\P^n$, also denoted $\omega_g$.
However, the affine open sets
$$U_i = \{ x \in \P^n : f_i(x) \ne 0\}$$
form an open cover $\U$ of $\P^n$.  Then $\omega_g$ is holomorphic on
$U_0\cap\cdots\cap U_n$, so it is a \v Cech cochain in
$C^n(\U,\Omega^n_{\P^n})$.  Further, since $\U$ has $n+1$ elements,
$\omega_g$ is a \v Cech cocycle and thus gives a class $[\omega_g]
\in H^n(\U,\Omega^n_{\P^n}) = H^n(\P^n,\Omega^n_{\P^n})$.  This is the
class $[\omega_g]$ mentioned in the introduction.
\medskip

	{\bf \S2. Residues on Toric Varieties.}  We will work with a
$n$-dimensional projective toric variety $X$ over the complex numbers
$\C$.  Thus $X$ is determined by a complete fan $\Sigma$ in $N_\R =
\R^n$.  As usual, $M$ denotes the dual lattice of $N = \Z^n$ and
$\Sigma(1)$ denotes the set of 1-dimensional cones in $\Sigma$.  Each
$\rho \in \Sigma(1)$ determines a divisor $D_\rho$ on $X$ and a
generator $n_\rho \in N\cap\rho$.  Standard references for toric
varieties are [D], [F] and [O].  As explained in [C], $X$ also has the
homogeneous coordinate ring $S = \C[x_\rho]$, which is graded by the
Chow group $A_{n-1}(X)$ so that a monomial $\Pi_\rho x_\rho^{a_\rho}$
has degree given by the class $[\sum_\rho a_\rho D_\rho] \in
A_{n-1}(X)$.  Given a class $\alpha \in A_{n-1}(X)$, we let $S_\alpha$
denote the graded piece of $S$ in degree $\alpha$, and we write
$\deg(f) = \alpha$ when $f \in S_\alpha$.

	Our strategy for defining toric residues is inspired by the
Trace Property of the introduction.  Thus we need, first, a trace map
$H^n(X,\Omega^n_X) \simeq \C$ and, second, a method that uses
polynomials $g \in S_\rho$ (for some $\rho \in A_{n-1}(X)$) to create
\v Cech cohomology classes $[\omega_g] \in H^n(X,\Omega^n_X)$.  Then
the toric residue will be easy to define.

	We begin with the trace map.  Since $X$ need not be smooth, we
can't use the usual sheaf on $n$-forms on $X$.  Instead, we use the
sheaf of Zariski $n$-forms on $X$, which by abuse of notation will be
written $\Omega^n_X$ (thus $\Omega_X^n = j_*\Omega^n_U$, where $j : U
\to X$ is the inclusion of the smooth part of $X$).  Since the toric
variety $X$ is Cohen--Macaulay with $\Omega_X^n$ as dualizing sheaf,
we have a trace map ${\rm Tr}_X : H^n(X,\Omega^n_X) \simeq \C$.  The duality
theory used here can be found in [O, \S3.2].

	Given an ample class $\beta \in A_{n-1}(X)$, there is a line
bundle $\ox(\beta)$ on $X$ and a canonical isomorphism $ S_\beta
\simeq H^0(X,\ox(\beta))$ (see [C, \S3]).  Regarding $f$ as a section
of $\ox(\beta)$, we can talk about what it means for $f$ to vanish at
point of $X$.  For the remainder of the paper, we make the following
assumption:
$$\hbox{$\beta \in A_{n-1}(X)$ is ample and $f_0,\dots,f_n \in
	S_\beta$ don't vanish simultaneously on $X$}\ .\leqno(3)$$

	Given $f_0,\dots,f_n$ as above, we set $U_i = \{x \in X :
f_i(x) \ne 0\}$.  Assumption 2.1 implies that the $U_i$ form an affine
open cover $\U$ of $X$.  As in \S1, we can use this open cover to
compute $H^n(X,\Omega_X^n)$ by \v Cech cohmology, so that every
section of $\Omega^n_X$ over $U_0\cap\cdots\cap U_n$ is a \v Cech
cocycle in $C^n(\U,\Omega^n_X)$.  Thus every $\omega \in
\Omega^n_X(U_0\cap\cdots\cap U_n)$ gives $[\omega] \in
H^n(X,\Omega^n_X)$.

	It remains to study sections of $\Omega_X^n$ over
$U_0\cap\cdots\cap U_n$.  We begin by constructing an analog of the
form (2).  Fix an integer basis $m_1,\dots,m_n$ for the lattice $M$.
Then, given a subset $I = \{\rho_1,\dots,\rho_n\} \subset \Sigma(1)$
consisting of $n$ elements, define
$$ \det(n_I) = \det(\langle m_i,n_{\rho_j}
	\rangle_{\scriptscriptstyle{1 \le i,j \le n}})\ . $$
Also set $dx_I = dx_{\rho_1}\wedge\cdots\wedge dx_{\rho_n}$ and
$\hat{x}_I = \Pi_{\rho \notin I} x_\rho$.  Note that $\det(n_I)$ and
$dx_I$ depend on how the $\rho \in I$ are ordered, while their product
$\det(n_I) dx_I$ does not.  Then we define the $n$-form $\Omega$ by
the formula
$$\Omega = \sum_{|I| = n} \det(n_I)\, \hat{x}_I\, dx_I\ ,\leqno(4)$$
where the sum is over all $n$-element subsets $I \subset
\Sigma(1)$.  This form is well-defined up to $\pm1$.

	Now consider the graded $S$-module $\widehat\Omega_S^n =
S\cdot\Omega$, where $\Omega$ is considered to have degree
$$\beta_0 = \textstyle{\sum_\rho} \deg(x_\rho) =
	[\textstyle{\sum_\rho} D_\rho] \in A_{n-1}(X)\ .$$
Thus $\widehat\Omega_S^n \simeq S(-\beta_0)$ as graded $S$-modules.
By [C, \S3], every graded $S$-module gives rise to a sheaf on $X$, and
by [BC, \S9], the sheaf associated to $\widehat\Omega^n_S$ is exactly
$\Omega_X^n$.  Furthermore, we can describe sections of $\Omega_X^n$
with prescribed poles as follows.

\proclaim Proposition 2.1. Let $\alpha \in A_{n-1}(X)$ be a Cartier
class, and let $Y \subset X$ be defined by the vanishing of $f \in
S_\alpha$.  Then
$$H^0(X,\Omega_X^n(Y)) = \Big\{ {g\,\Omega \over f} : g \in
	S_{\alpha-\beta_0}\Big\} \simeq S_{\alpha - \beta_0}\ .$$

\prf This follows from Proposition 9.7 of [BC].  (Although [BC]
assumes that X is simplicial, the results of \S8 and \S9 of [BC] apply
to all complete toric varieties.) \qed
\medskip

	If we apply this proposition to $f = f_0\cdots f_n \in
S_{(n+1)\beta}$, we get an $n$-form
$$\omega_g = {g\,\Omega \over f_0\cdots f_n} \in
	\Omega_X^n(U_0\cap\cdots\cap U_n)$$
for all $g \in S_{(n+1)\beta - \beta_0}$.  To simplify notation, we set
$$ \rho = (n+1)\beta - \beta_0\ .$$
Hence there are classes $[\omega_g] \in H^n(X,\Omega_X^n)$ for all $g
\in S_\rho$.

	Now we can finally define the toric residue.

\proclaim Definition 2.2. If $f_0,\dots,f_n \in S_\beta$ satisfy (3)
and $g \in S_\rho$, the {\bf toric residue} is
$$\res(\omega_g) = {\rm Tr}_X([\omega_g])\ .$$

	The first properties of toric residues are easy to prove.

\proclaim Proposition 2.3.
\vskip0pt
\item{(i)} $\res(\omega_g)$ is $\C$-linear in $g$ and antisymmetric in
$f_0,\dots,f_n$.
\item{(ii)} $\res(\omega_g) = 0$ whenever $g \in \langle
f_0,\dots,f_n\rangle_\rho$.
\vskip0pt

\prf The linearity of the toric residue is obvious.  The isomorphism
$H^n(\U,\Omega^n_{\P^n}) \simeq H^n(\P^n,\Omega^n_{\P^n})$ comes from
the differential in the \v Cech complex and hence depends
asymmetrically on how we order the open sets $U_i = \{f_i \ne 0\}$ of
${\cal U}$.  Thus the toric residue is asymmetric in the $f_i$.  The
second part of the proposition follows easily by considering the \v
Cech coboundary $\delta : C^{n-1}(\U,\Omega^n_X) \to
C^n(\U,\Omega_X^n)$. We omit the details.\qed
\medskip

	As a corollary, we see that the toric residue induces a map
$$\res : S_\rho/\langle f_0,\dots,f_n\rangle_\rho \longrightarrow \C
	\ .$$
In the next section, we will see that the quotient on the right is
one-dimensional, and in \S5 we will show that the above map is in
fact an isomorphism.

	We should also mention that in the forthcoming paper [CCD], it
will be shown that when $X$ is smooth, the toric residue
$\res(\omega_g)$ can be expressed as a sum of local residues on $X$.
Namely, given $f_0,\dots,f_n \in S_\beta$ which satisfy (3), let
$D_i$ be the divisor $f_i = 0$.  Then for each $0 \le k \le n$, the
intersection $D_0\cap \cdots\cap\widehat{D_k}\cap\cdots \cap D_n$ is
finite since it is contained in the affine variety $X - D_k$.
This means for $g \in S_\rho$, the meromorphic form
$$\omega_g = {(g/f_k)\,\Omega \over f_0\cdots\widehat{f_k}\cdots
f_n}$$
has $D_0\cup\cdots\cup\widehat{D_k}\cup\cdots\cup D_n$ as local polar
divisor near $x \in D_0\cap\cdots\cap \widehat{D_k}\cap\cdots\cap
D_n$.  Thus, for each such $x$, we can define the Grothendieck local
residue symbol
$$ \res_x\left( { (g/f_k)\,\Omega \over
	f_0\cdots \widehat{f_k} \cdots f_n} \right)\ .$$
Then the following equality is a special case of the results of [CCD]:
$$\res(\omega_g) = (-1)^k \sum_{x \in D_0\cap\cdots\cap \widehat{D_k}
	\cap \cdots\cap D_n} \res_x\left( { (g/f_k)\,\Omega \over
	f_0\cdots \widehat{f_k} \cdots f_n} \right)\ .$$
\smallskip

	{\bf \S3. Some Commutative Algebra.} The basic commutative
algebra for our situation is given in [Bat, Theorem 2.10 and
Proposition 9.4].  In this section, we recast Batyrev's results in
terms of the ring $S$ and supply some of the details.

	We first study the subring of $S$ determined by an ample class
$\beta \in A_{n-1}(X)$.

\proclaim Proposition 3.1. If $\beta \in A_{n-1}(X)$ is ample, then
the ring $S_{*\beta}= \oplus_{k=0}^\infty S_{k\beta}$ is
Cohen--Macaulay of dimension $n+1$, with canonical module given by
$\omega_{S_{*\beta}} = \oplus_{k=0}^\infty S_{k\beta-\beta_0}$.

\prf Write $\beta$ as $[\sum_\rho a_\rho D_\rho]$.  Then $\Delta = \{
m \in M_\R : \langle m,n_\rho\rangle \ge -a_\rho\}$ is a
$n$-dimensional convex polyhedron since $\beta$ is ample.  Let
$\check\sigma \subset \R\times M_\R$ be the cone over $\{1\}\times
\Delta$.  The dual of $\check\sigma$ is a strongly convex rational
polyhedral cone $\sigma \subset \R\times N_\R$.  Then $\sigma$
determines a $(n+1)$-dimensional affine toric variety with coordinate
ring
$$\eqalign{S_\Delta = \C[\check\sigma\cap(\Z\times M)] &= \C[t_0^k t^m
	: (k,m) \in \check\sigma\cap (\Z\times M)]\cr
	&= \C[t_0^k t^m : (1,m/k) \in \{1\}\times\Delta]\cr
	&= \C[t_0^k t^m : m \in k\Delta]\ .}$$
Hence $S_\Delta$ is Cohen--Macaulay by [D, Theorem 3.4] or [Ho].  But
$S_\Delta \simeq S_{*\beta}$ follows from the proof of [BC, Theorem
11.5], so that $S_{*\beta}$ is Cohen--Macaulay.

	 We next determine the canonical module of $S_{*\beta}$ (see
[BH] for more background on this topic).  By [D, 4.6], the canonical
module of the semigroup ring $S_\Delta$ is a certain ideal
$I^{(1)}_\Delta \subset S_\Delta$ (this is the notation of [Bat,
\S2]).  Then the proof of [BC, Theorem 11.8] shows that under the
isomorphism $S_\Delta \simeq S_{*\beta}$, the ideal $I^{(1)}_\Delta
\subset S_\Delta$ maps to $\oplus_{k=0}^\infty \langle \Pi_\rho x_\rho
\rangle_{k\beta} \subset S_{*\beta}$.  This is isomorphic to
$\oplus_{k=0}^\infty S_{k\beta-\beta_0}$ since $\Pi_\rho x_\rho$ has
degree $\beta_0$.

	Alternatively, notice that $S_{*\beta} = \oplus_{k=0}^\infty
S_{k\beta}$ has a natural grading and that $X = {\rm
Proj}(S_{*\beta})$ (see [Bat, \S2]).  Then the canonical module is
$\oplus_{k=0}^\infty H^0(X,\Omega^n_X(k\beta))$ by [GW, 5.1.8].
However, by Proposition 2.1, we can identify
$H^0(X,\Omega^n_X(k\beta))$ with $S_{k\beta-\beta_0}$, which shows
that $\oplus_{k=0}^\infty S_{k\beta-\beta_0}$ is the canonical
module. \qed
\medskip

	For a general Cohen--Macaulay variety $X$ and ample
$\beta$, the ring $\oplus_{k=0}^\infty H^0(X,\ox(k\beta))$ need not be
Cohen--Macaulay, although some Veronese subring will be (see [GW,
5.1.11]).

	We next bring $f_0,\dots,f_n \in S_\beta$ into the picture.

\proclaim Proposition 3.2. If $f_0,\dots,f_n \in S_\beta$ satisfy (3),
then:
\vskip0pt
\item{(i)} $f_0,\dots,f_n$ are a regular sequence in $S_{*\beta} =
\oplus_{k=0}^\infty S_{k\beta}$.
\item{(ii)} $R = S_{*\beta}/\langle f_0,\dots,f_n\rangle$ is a
zero-dimensional Cohen--Macaulay ring.
\item{(iii)} The canonical module of $R$ is $\omega_R =
(\omega_{S_{*\beta}}/\langle
f_0,\dots,f_n\rangle\omega_{S_{*\beta}})[n+1]$ (where the $[n+1]$
indicates a shift in grading).
\vskip0pt

\prf Since $f_0,\dots,f_n$ define the empty subvariety of $X$, $R =
S_{*\beta}/\langle f_0,\dots,f_n\rangle$ has dimension zero.  Then (i)
and (ii) follow from [BH, Theorems 2.1.2 and 2.1.3], and (iii) follows
from [BH, Corollary 3.6.14].  \qed
\medskip

	We can now determine $S_\rho/\langle
f_0,\dots,f_n\rangle_\rho$, where as usual $\rho = (n+1)\beta -\beta_0$.

\proclaim Proposition 3.3.  There is a natural isomorphism
$S_\rho/\langle f_0,\dots,f_n\rangle_\rho \simeq \C$.

\prf This is follows from local duality for graded Cohen--Macaulay
rings (see [BH, 3.6]).  In the zero-dimensional case, local duality
is a natural isomorphism of graded $R$-modules
$$\omega_R \simeq {\rm Hom}_\C(R,\C)\ .$$
In particular, $(\omega_R)_0 \simeq {\rm Hom}_\C(R_0,\C) \simeq \C$
since $R_0 = \C$.  By Propositions 3.1 and 3.2, we have
$$ (\omega_R)_0 = (\omega_{S_{*\beta}}/\langle f_0,\dots,f_n\rangle
	\omega_{S_{*\beta}})_{n+1} = S_{(n+1)\beta-\beta_0}/\langle
	f_0,\dots,f_n \rangle_{(n+1)\beta-\beta_0}\ ,$$
and the corollary follows. \qed
\medskip

	Note that when $\beta = \beta_0$ (so $X$ is a Fano toric
variety), we have $\omega_R \simeq R[n]$, so that $S_{*\beta_0}$ and
$R$ are Gorenstein (see [BH, 3.6.11]).  This case is of interest in
mirror symmetry.
\medskip

	{\bf \S4. Toric Jacobians.} This section will define a
``Jacobian'' of $f_0,\dots,f_n \in S_\beta$ which is closely related
to the Jacobian $\det(\partial f_i/\partial x_j)$ when $X = \P^n$.
Here is our main result.

\proclaim Proposition 4.1. If $f_0,\dots,f_n \in S_\alpha$, then there
is $J \in S_{(n+1)\alpha - \beta_0}$ such that
$$\textstyle{\sum_{i=0}^n} (-1)^i f_i\, df_0\wedge \cdots \wedge
	\widehat{df_i}\wedge \cdots \wedge df_n = J\,\Omega\ .$$
Furthermore, if $I = \{\rho_1,\dots,\rho_n\} \subset \Sigma(1)$ and
$n_{\rho_1},\dots,n_{\rho_n}$ are linearly independent, then
$$ J = { \det\pmatrix{f_0 & \cdots & f_n\cr
	{\partial f_0/\partial x_{\rho_1}} & \cdots &
	{\partial f_n / \partial x_{\rho_1}}\cr
	\vdots & & \vdots \cr
	{\partial f_0 / \partial x_{\rho_n}} & \cdots &
	{\partial f_n / \partial x_{\rho_n}}\cr}
	\over \det(n_I)\, \hat{x}_I}\ .\leqno(5)$$

\prf Note that we assume nothing about the ampleness of $\alpha$ or
the vanishing of the $f_i$'s.  This proof was suggested by Eduardo
Cattani and Alicia Dickenstein.  We can assume $f_0 \ne 0$, so that
$$\textstyle{\sum_{i=0}^n} (-1)^i f_i\, df_0\wedge \cdots \wedge
	\widehat{df_i}\wedge \cdots \wedge df_n = f_0^{n+1}
	d(f_1/f_0)\wedge\cdots\wedge d(f_n/f_0)\ .$$
Let $m_1,\dots,m_n$ be the basis of $M$ used in the definition of
$\Omega$, and set
$$t_i = \Pi_\rho x^{\langle m_i,n_\rho\rangle}_\rho\ .$$
Then $t_1,\dots,t_n$ are coordinates for the torus $T \subset X$ (see
[BC, \S9]), and each $f_i/f_0$ is a rational function of the $t_i$
since $f_i$ and $f_0$ have the same degree.  Hence
$$d(f_1/f_0)\wedge\cdots\wedge d(f_n/f_0) = J(t_1,\dots,t_n)\,
	dt_1\wedge\cdots\wedge dt_n$$
for some rational function $J(t_1,\dots,t_n)$.  However, the proof of
Proposition 9.5 in [BC] shows that
$$\Omega = \Pi_\rho x_\rho\, {dt_1\over t_1}\wedge\cdots\wedge
	{dt_n\over t_n}\ ,$$
and it follows easily from the above equations that
$$\textstyle{\sum_{i=0}^n} (-1)^i f_i\, df_0\wedge \cdots \wedge
	\widehat{df_i}\wedge \cdots \wedge df_n = J\,\Omega\ ,$$
where
$$ J = f_0^{n+1}\,J(t_1,\dots,t_n)\,{t_1\cdots t_n\over\Pi_\rho
	x_\rho}\ .$$
This equation shows that $J$ has degree $(n+1)\alpha - \beta_0$ as a
rational function of the $x_\rho$'s.

	If $I \subset \Sigma(1)$ has $|I| = n$, we let $J(f_I)$ denote
the determinant in the numerator of (5).  Then one easily computes
that
$$\textstyle{\sum_{i=0}^n} (-1)^i f_i\, df_0\wedge \cdots \wedge
	\widehat{df_i}\wedge \cdots \wedge df_n =
	\textstyle{\sum_{|I|=n}} J(f_I)\, dx_I\ .$$
Since the right hand side equals $J\Omega$, it follows that
$$ J\,\det(n_I)\,\hat{x}_I = J(f_I)$$
for all $I$.  This gives the desired formula (5) for $J$.

	It remains to show that $J$ is a polynomial in the
$x_\rho$'s.  If we write $J$ as a quotient of relatively prime
polynomials in $S$, then (5) shows that the denominator divides
$\hat{x}_I$ for every $I$ with $\det(n_I) \ne 0$.  Since the fan of
$X$ is complete, every $\rho \in \Sigma(1)$ is in some such $I$
($\rho$ lies in a $n$-dimensional cone $\sigma$, and $I$ can be
chosen to be an appropriate subset of $\sigma(1)$ containing $\rho$).
It follows that the $\hat{x}_I$'s are relatively prime, which forces
the denominator of $J$ to be a constant.  Since we've already seen
that $J$ has degree $(n+1)\alpha - \beta_0$, it follows immediately
that $J\in S_{(n+1)\alpha - \beta_0}$, and the proposition is
proved.\qed
\medskip

	In light of this proposition, we make the following
definition.

\proclaim Definition 4.2.  Given $f_0,\dots,f_n \in S_\alpha$, the
polynomial $J \in S_{(n+1)\alpha-\beta_0}$ satisfying
$$\textstyle{\sum_{i=0}^n} (-1)^i f_i\, df_0\wedge \cdots \wedge
	\widehat{df_i}\wedge \cdots \wedge df_n = J\,\Omega$$
is called the {\bf toric Jacobian} of $f_0,\dots,f_n$.

	For example, suppose $X = \P^n$ and $f_0,\dots,f_n$ are
homogeneous of degree $d$.  If $I$ is given by $x_1,\dots,x_n$, then
applying the Euler formula $f_j = (1/d)\sum_{i=0}^n x_i \partial
f_j/\partial x_i$ to (5) shows that the toric Jacobian is given by
$$ \eqalign{J &=  {\det\pmatrix{f_0 & \cdots & f_n\cr
	{\partial f_0/\partial x_{1}} & \cdots &
	{\partial f_n / \partial x_{1}}\cr
	\vdots & & \vdots \cr
	{\partial f_0 / \partial x_{n}} & \cdots &
	{\partial f_n / \partial x_{n}}\cr}
	\over x_0} =
	{\det\pmatrix{{1\over d}x_0\partial f_0/\partial x_0 & \cdots &
	{1\over d}x_0 \partial f_n/\partial x_0\cr
	{\partial f_0/\partial x_{1}} & \cdots &
	{\partial f_n / \partial x_{1}}\cr
	\vdots & & \vdots \cr
	{\partial f_0 / \partial x_{n}} & \cdots &
	{\partial f_n / \partial x_{n}}\cr}
	\over x_0}\cr
	&= {1\over d} \det(\partial f_i/\partial x_j)\ .\cr}$$
In [PS, 12.10], it was assumed that $J = \det(\partial f_i/\partial
x_j)$, which caused the residue formula given there to have an extra
factor of $d$.

	The toric Jacobian is also related to the hyperdeterminant of
a $m\times(m+p-1)\times p$ matrix $A = (a_{ijk})$, as described in
[GKZ, \S3 of Chapter 14].  From $A$, we get $m+p-1$ bilinear forms
$f_j = \sum_{ik} a_{ijk} x_i y_k$ in variables
$x_1,\dots,x_m,y_1,\dots,y_p$.  Thus $X = \P^{m-1}\times\P^{p-1}$ and
$f_j \in S_{1,1}$, where $S = \C[x_i;y_k]$ has the usual bigrading.
The toric Jacobian $J$ of $f_1,\dots,f_{m+p-1}$ has degree
$(m+p-1)(1,1) - (m,p) = (p-1,m-1)$.  In [GKZ, Chapter 14], $J$ appears
in equation (3.20), and in Theorem 3.19, the coefficients of $J$ are
used to compute the hyperdeterminant of $A$.  Also, Proposition 3.21
gives an interesting combinatorial interpretation of the coefficients
of $J$.

		As our final example, let $f =
x^2z^2+x^2w^2+y^2z^2+y^2w^2+\lambda xyzw \in \C[x,y;z,w]$.  Then $f$
has degree $(2,2)$, so that the toric Jacobian of $f, xf_x, zf_z$
has degree $3(2,2) - (2,2) = (4,4)$.  Since $\Omega =
(xdy-ydx)\wedge(zdw-wdz)$ for this case, we know the $\det(n_I)$'s,
and one computes that the toric Jacobian equals
$$4(\lambda x^4z^2w^2+4x^3yzw^3+4x^3yz^3w+\lambda x^2y^2z^4 +
	\lambda x^2y^2w^4 + 4xy^3z^3w + 4xy^3zw^3 + \lambda y^4z^2w^2)
	\ .$$
Other examples are equally easy to compute.
\medskip

	{\bf \S5. Uniqueness of Toric Residues.} Putting together what
we proved in the previous sections, we can now state the main theorem
of this paper.

\proclaim Theorem 5.1. Let $X$ be a complete toric variety, and let
$\beta \in A_{n-1}(X)$ be ample.  If $f_0,\dots,f_n \in S_\beta$ don't
vanish simultaneously on $X$, then
\vskip0pt
\item{(i)} If $\rho = (n+1)\beta -\beta_0$, the toric residue map
$\res : S_\rho/\langle f_0,\dots,f_n\rangle_\rho \to \C$ from \S2 is
an isomorphism.
\item{(ii)} If $J \in S_\rho$ is the toric Jacobian of $f_0,\dots,f_n$
from \S4, then
$$ \res(\omega_J) = n!\,{\rm vol}(\Delta) = \deg(F)\ ,$$
where ${\rm vol}(\Delta)$ is the normalized volume of the convex
polyhedron $\Delta \subset M_\R$ associated to $\beta$ (see the proof
of Proposition 3.1) and $F : X \to \P^n$ is the map defined by $F(x) =
(f_0(x),\dots,f_n(x))$.
\vskip0pt

\prf We know from Proposition 3.3 that $S_\rho/\langle
f_0,\dots,f_n\rangle_\rho$ has dimension one.  Hence (i) is an
immediate consequence of (ii).  To prove (ii), note that $F =
(f_0,\dots,f_n) : X \to \P^n$ comes from $n+1$ sections of an ample
line bundle.  Since the sections never vanish simultaneously and ${\rm
dim}(X) = n$, $F$ is defined everywhere and is finite and surjective.

	We now proceed as in [PS, 12.10].  By the definition of the
toric Jacobian, we have
$$\eqalign{J\,\Omega &= \textstyle{\sum_{i=0}^n} (-1)^i f_i\,
	df_0\wedge \cdots \wedge \widehat{df_i}\wedge \cdots \wedge
	df_n\cr &=
	F^*(\textstyle{\sum_{i=0}^n} (-1)^i x_i\, dx_0\wedge \cdots
	\wedge \widehat{dx_i}\wedge \cdots \wedge dx_n) =
	F^*(\Omega_{\P^n})\ .\cr}$$
Thus
$$\omega_J = {J\,\Omega \over f_0\cdots f_n} = F^*\Big({\Omega_{\P^n}
	\over x_0\cdots x_n}\Big)\ .$$
Denote the $n$-form in parentheses by $\omega_1$.  Then the map $F^* :
H^n(\P^n,\Omega^n_{\P^n}) \to H^n(X,\Omega_X^n)$ satisfies
$F^*([\omega_1]) = [\omega_J]$ (this is easy to see using \v Cech
cohomology).  Since $F$ is finite and surjective, standard properties
of the trace map imply that
$$\res(\omega_J) = {\rm Tr}_X([\omega_J]) = {\rm
	Tr}_X(F^*([\omega_1])) = \deg(F)\,{\rm Tr}_{\P^n}([\omega_1])
	= \deg(F)$$
(see [H, Chapter III]).

	To complete the proof, we need to show that $\deg(F) = n!\,{\rm
vol}(\Delta)$.  If $D$ is the divisor of a section of the line bundle
$\ox(\beta)$, then it is well-known that
$$ D^n = n!\,{\rm vol}(\Delta)$$
(see [O, Proposition 2.10]).  Note that when $X$ is singular, we use
the intersection number $D^n$ as defined in [K, Chapter I] or [O,
\S2.2].  Since $F : X \to \P^n$ is finite and $F^*({\cal O}_{\P^n}(1))
\simeq \ox(\beta)$, it follows from [K, Proposition 6 of Chapter I,
\S2] that
$$D^n = \deg(F)\,H^n = \deg(F)\ ,$$
where $H \subset \P^n$ is a hyperplane.  Thus $\deg(F) = n!\,{\rm
vol}(V)$, and we are done.\qed
\medskip

	Hence toric residues have both the Quotient and Trace
Properties mentioned in the introduction.  Note also that the toric
residue is uniquely characterized by these properties.

	When $X = \P^n$ and $f_0,\dots,f_n$ have degree $d$, we saw in
\S4 that the toric Jacobian is $J = (1/d)\det(\partial f_i/\partial
x_j)$.  Since $F = (f_0,\dots,f_n)$ has degree $d^n$, if we set $g =
\det(\partial f_i/\partial x_j)$, then
$$\res(\omega_g) = d\,\res(\omega_J) = d\,\deg(F) = d^{n+1}\ .$$
However, for this choice of $g$, the Grothendieck residue symbol (1)
is the local intersection number of the divisors in $\C^{n+1}$ defined
by $f_i = 0$, which is also $d^{n+1}$.  Thus
$$\res(\omega_g) = \res_{0}\left( {g\, dx_0\wedge \cdots\wedge dx_n
	\over f_0\cdots f_n} \right)\ ,$$
and from here, one easily sees that equality holds for all $g \in S_\rho$
(see [PS, 12.10]).

	For another application of our theory, suppose that $f \in
S_\beta$, where $\beta$ is ample and $f$ is nondegenerate in the sense
of [BC, Definition 4.13].  In this situation, $f$ defines a
hypersurface $Y \subset X$, and the ideal
$$J_0(f) = \langle x_\rho \partial f/\partial x_\rho\rangle \subset
	S$$
is closely related to the mixed Hodge structure of the affine
hypersurface $Y\cap T$, where $T \subset X$ is the torus of $X$ (see
[Bat, \S9] or [BC, \S11]).  Then we get the following proposition.

\proclaim Proposition 5.3. Given $f \in S_\beta$, suppose that
$n_{\rho_1},\dots,n_{\rho_n}$ are linearly independent.  Then
$J_0(f) = \langle f, x_{\rho_1}\partial f/\partial
x_{\rho_1},\dots,x_{\rho_n}\partial f/\partial x_{\rho_n} \rangle$.
Furthermore,
$$\eqalign{f\ \hbox{is nondegenerate}\ &\iff\
	x_\rho \partial f/\partial x_\rho\ \hbox{don't vanish
	simultaneously on $X$}\cr &\iff\
	f,\ x_{\rho_i} \partial f/\partial x_{\rho_i}\ \hbox{don't
	vanish simultaneously on $X$\ .}\cr}$$
Finally, if $\beta$ is ample and $f$ is nondegenerate, then
$S_\rho/J_0(f)_\rho \simeq \C$ ($\rho = (n+1)\beta -\beta_0$), and the
toric Jacobian of $f,x_{\rho_1}\partial f/\partial x_{\rho_1}
\dots,x_{\rho_n}\partial f/\partial x_{\rho_n}$ represents a nonzero
element of $S_\rho/J_0(f)_\rho$.

\prf Our proof will use Euler formulas.  Recall from [BC, \S3]
that $\theta = \sum_\rho b_\rho x_\rho\partial/\partial x_\rho$ is an
Euler vector field provided $\sum_\rho b_\rho n_\rho = 0$.  Further,
such a $\theta$ determines a constant $\theta(\beta)$ such that the
Euler formula
$$\theta(\beta)\,f = \textstyle{\sum_\rho} b_\rho x_\rho\partial f
	/\partial x_\rho$$
holds for all $f \in S_\beta$.

	Now, given any $\rho$, there is a relation $n_\rho +
\sum_{i=1}^n b_i n_{\rho_i} = 0$ since the $n_{\rho_i}$ are a basis
over $\Q$.  The resulting Euler formula shows that $x_\rho\partial
f/\partial x_\rho \in \langle f,x_{\rho_i} \partial f/\partial
x_{\rho_i}\rangle$.  We also have $f \in J_0(f)$ (see the proof of
[BC, Lemma 10.5]), and it follows that $f$ and the $x_{\rho_i}\partial
f/\partial x_{\rho_i}$ generate $J_0(f)$.

	The argument of [BC, Proposition 3.5] adapts easily to show
that $f$ is nondegenerate if and only if $x_\rho
\partial f/\partial x_\rho$ don't vanish simultaneously.  This proves
the first equivalence of the proposition, and second equivalence is
now trivial.

	Finally, when $\beta$ is ample and $f \in S_\beta$ is
nondegenerate, we can apply Theorem 5.1 to $f,x_{\rho_1}\partial
f/\partial x_{\rho_1},\dots,x_{\rho_n}\partial f/\partial
x_{\rho_n}$.  This gives the final part of the proposition.\qed
\medskip

	An example of this can be found in \S4, where $X =
\P^1\times\P^1$ and $f = x^2z^2+x^2w^2+y^2z^2+y^2w^2+\lambda xyzw \in
S = \C[x,y;z,w]$.  One can check that $f$ is nondegenerate provided
$\lambda \ne 0, \pm4$, so that the toric Jacobian of $f$, $x\partial
f/\partial x$ and $z\partial f/\partial z$ (displayed at the end of
\S4) is a nonzero element of $S_{4,4}/J_0(f)_{4,4} \simeq \C$.

	In general, the isomorphism $S_\rho/J_0(f)_\rho \simeq \C$ was
previously known to Batyrev and is related to the cup-product pairing
$$ H^{n-1}(Y\cap T,\C) \otimes H^{n-1}_c(Y\cap T,\C) \longrightarrow
	\C$$
(see [Bat, \S9] or [BC, \S11]).  If we assume in addition that $X$ is
simplicial, then the ideal quotient
$$ J_1(f) = J_0(f) : \textstyle{\Pi_\rho} x_\rho$$
also plays a role: it is related to the Hodge structure on the
primitive cohomology of the hypersurface $Y \subset X$ defined by $f$
(still assuming that  $f \in S_\beta$ is nondegenerate and $\beta$
ample).  In this case, one can show that natural map
$$ S_{(n+1)\beta - 2\beta_0}/J_1(f)_{(n+1)\beta-2\beta_0}
	\mapname{\Pi_\rho x_\rho} S_\rho/J_0(f)_\rho$$
is an isomorphism, so that we have an isomorphism
$$ S_{(n+1)\beta - 2\beta_0}/J_1(f)_{(n+1)\beta-2\beta_0} \simeq \C
	\ .\leqno(6)$$
This isomorphism is closely related to cup product on the primitive
cohomology of $Y$.

	Is there an explicit formula for a polynomial that gives a
nonzero element in the quotient (6)?  In the case of $X = \P^n$, the
Hessian $\partial^2 f/\partial x_i\partial x_j$ is such an element.
It would be very interesting to have a toric generalization of the
Hessian.
\medskip

	{\bf \S6. Integral Representations.} In this section, we will
explore several ways of representing the toric residue as an
integral.  In order to do this, we will assume that $X$ is a
projective simplicial toric variety.  Thus $X$ is an orbifold (or
$V$-manifold), so that by [Bai], we have a Dolbeault isomorphism
$H_{\dbar}^{n,n}(X) \simeq H^n(X,\Omega^n_X)$.

	Our first set of formulas use the Dolbeault isomorphism to
express $\res(\omega_g)$ as an integral.  We begin with $f_0,\dots,f_n
\in S_\beta$ which satisfy (3).  As in \S2, we get the open covering
$\U$ of $X$ and the \v Cech cocyle $\omega_g \in C^n(\U,\Omega_X^n)$
for $g \in S_\rho$.  Now let
$$\rho_i = {|f_i|^2 \over |f_0|^2 + \cdots + |f_n|^2}$$
for $0 \le i\le n$.  Since the $f_i$ all have the same degree and
don't vanish simultaneously, the $\rho_i$ are $C^\infty$-functions on
$X$ which sum to $1$.

\proclaim Proposition 6.1. The class of the $C^\infty$ $(n,n)$-form
$$ \eta_g = (-1)^{n(n+1)/2} n!\,{g\,\textstyle{\sum_{i=0}^n} (-1)^i
	\overline{f}_i\, d\overline{f}_0\wedge\cdots\wedge
	\widehat{d\overline{f}_i}\wedge\cdots\wedge d\overline{f}_n
	\wedge\Omega \over (|f_0|^2 + \cdots + |f_n|^2)^{n+1}}$$
maps to $[\omega_g]$ under the Dolbeault isomorphism
$H_{\dbar}^{n,n}(X) \simeq H^n(X,\Omega^n_X)$.

\noindent {\bf Remark.} This proposition gives a more precise version
of a formula in [GH, Chapter 5].  There, the factor
$(-1)^{n(n+1)/2}n!$ appears as a constant $C_n$.
\medskip

\prf We first show that $\eta_g$ lives on $X$.  Since $X$ is
simplicial and complete, we can write $X$ as a geometric quotient
$(\C^{\Sigma(1)} - Z)/G$, where $Z \subset \C^{\Sigma(1)}$ is a proper
closed subvariety and $G = {\rm Hom}_\Z(A_{n-1}(X),\C^*)$ (see [C,
\S2]).  This is also a quotient in the real analytic category, so that
the form $\eta_g$ on $\C^{\Sigma(1)} - Z$ descends to $X$ if and only
if it is invariant under $G$ and is annihilated by all real vectors
fields to $G$.  From [BC, \S3], we know that the complex Lie algebra
of $G$, denoted ${\rm Lie}(G)$, consists of the Euler vector fields
$\theta = \sum_\rho b_\rho x_\rho \partial /\partial x_\rho$ described
in the proof of Proposition 5.2.

	The form $\eta_g$ is clearly invariant under $G$ since all
$f_i$ have the same degree.  To see that it is annihilated by all real
vector fields to $G$, it suffices to show that
$$\theta\ip\eta_g = \overline{\theta}\ip\eta_g = 0$$
for all $\theta \in {\rm Lie}(G)$.  However, Proposition 4.1 implies
that $\eta_g$ is a $C^\infty$ multiple of $\overline{\Omega}
\wedge\Omega$.  Thus we need to prove that
$$\theta\ip\Omega = \overline{\theta}\ip\overline{\Omega} =
\theta\ip\overline{\Omega} = \overline{\theta}\ip\Omega = 0\ .$$
The last two of these are trivially zero, and the vanishing of the
second follows from the vanishing of the first.  Hence it remains to
show that $\theta\ip\Omega = 0$ for any $\theta \in {\rm Lie}(G)$.
This will follow immediately from the following observation of Batyrev.

\proclaim Lemma 6.2. Let $\theta_1,\dots,\theta_r$ be an ordered
basis of ${\rm Lie}(G)$, and let $d{\bf x} = \wedge_\rho dx_\rho$ for
some ordering of the $x_\rho$'s.  Then there is a nonzero constant $c$
such that
$$(\theta_1\wedge\cdots\wedge\theta_r)\ip d{\bf x} = c\,\Omega$$
(where $\ip$ denotes interior multiplication).

\prf  Consider the map $\oplus_\rho \C x_\rho \partial/\partial x_\rho
\to N_\C$ which sends $x_\rho \partial/\partial x_\rho$ to $n_\rho$.
This map is onto since $X$ is complete, and the kernel consists of the
Euler vector fields, which as above form the Lie algebra ${\rm
Lie}(G)$.  Thus we have an exact sequence
$$0 \longrightarrow {\rm Lie}(G) \longrightarrow
	\textstyle{\oplus_\rho} \C\, x_\rho \partial/\partial x_\rho
	\longrightarrow N_\C \longrightarrow 0\ .\leqno(7)$$

	We first compute $(\theta_1\wedge\cdots\wedge\theta_r)\ip
d{\bf x}$.  Write $\theta_i = \sum_\rho b^i_\rho x_\rho
\partial/\partial x_\rho$ for $1 \le i \le r$, and if $I \subset
\Sigma(1)$ has cardinality $n$, let
$$ \det(\hat{b}_I) = \det(b^i_\rho : 1 \le i \le r,\ \rho \notin
	I)\ .$$
The matrix $(b^i_\rho)$ is square since $r = |\Sigma(1)| - n$ by
(7).  Then one computes that
$$(\theta_1\wedge\cdots\wedge\theta_r)\ip d{\bf x} =
	\textstyle{\sum_{|I|=n}} (-1)^I \det(\hat{b}_I)\, \hat{x}_I\,
	dx_I\ ,$$
where $\hat{x}_I\,dx_I$ is as in \S2 and $(-1)^I$ is the sign of the
permutation of $\Sigma(1)$ which puts the $\rho \in I$ at the end but
otherwise preserves their order (see [S, p.~21]).

	To relate this to $\Omega$, we will compute the the
determinant (as in [GKZ, Appendix A]) of the exact sequence (7) with
respect to the following ordered bases.  First,
$\theta_1,\dots,\theta_r$ give an ordered basis of ${\rm Lie}(G)$, and
the ordering of the $x_\rho$ gives an obvious ordered basis of the
middle term of (7).  Finally, recall that in \S2 we picked an ordered
basis $m_1,\dots,m_n$ of $M$.  So we let $m_1^*,\dots,m_n^*$ be the
dual basis of $N$.  Then, by [GKZ, Appendix A], the based exact
sequence (7) has a determinant $c \in \C^*$.  We claim that
$$ c\, \det(n_I) = (-1)^I\det(\hat b_I)\ .\leqno(8)$$
for all $I \subset \Sigma(1)$ with $|I| = n$.  Comparing the formula
(4) for $\Omega$ to the above formula for
$(\theta_1\wedge\cdots\wedge\theta_r)\ip d{\bf x}$, the lemma will
follow immediately once (8) is proved.

	Let $I = \{\rho_1,\dots,\rho_n\}$.  First note that (8) holds
when $n_{\rho_1},\dots,n_{\rho_n}$ are linearly dependent, since
$\det(n_I) = 0$, and $\det(\hat{b}_I) = 0$ is also true because any
relation $\sum_{i=0}^n b_i n_{\rho_i} = 0$ must be a consequence of
the relations $\sum_\rho b^i_\rho n_\rho = 0$.  On the other hand, if
$n_{\rho_1},\dots,n_{\rho_n}$ are linearly independent, consider the
ordered basis of $\oplus \C x_\rho \partial/\partial x_\rho$ where the
$x_\rho \partial/\partial x_\rho$ for $\rho \in I$ all appear at the
end (but otherwise we preserve their order).  Using this basis and the
bases $\theta_i$ and $m^*_i$ for the other terms of (7), Proposition
11 of [GKZ, Appendix A] expresses the determinant as a quotient of two
determinants, where the numerator is the $r\times r$ determinant
obtained using the $\theta_i$ and the $x_\rho \partial/\partial
x_\rho$ for $\rho \notin I$, and the denominator is the $n\times n$
determinant using $x_\rho \partial/\partial x_\rho$ for $\rho \in I$
and the $m^*_i$.  Hence the determinant of this based exact sequence
is
$$ c' = {\det(\hat{b}_I) \over \det(n_I)}\ .$$
However, $c$ and $c'$ differ by the change of basis determinant for
the middle term of (7) (see Proposition 9 of [GKZ, Appendix A]), so
that $c = (-1)^I c'$, and (8) is proved.\qed
\medskip

	To continue the proof of Proposition 6.1, we next recall the
Dolbeault isomorphism.  As in [GH, p.~45], we have the exact sequences
$$0 \longrightarrow {\cal Z}^{n,n-p-1} \longrightarrow {\cal
	A}^{n,n-p-1} \mapname{\dbar} {\cal Z}^{n,n-p} \longrightarrow
	0 \ ,$$
and the Dolbeault isomorphism $H_{\dbar}^{n,n}(X) \simeq
H^n(X,\Omega^n_X)$ is obtained by composing the coboundary maps
$$\delta_p : H^p(X,{\cal Z}^{n,n-p}) \longrightarrow H^{p+1}(X,{\cal
	Z}^{n,n-p-1})\ .$$
To prove the proposition, it suffices to show that
$\delta_{n-1}\circ\cdots\circ\delta_0(\eta_g) = [\omega_g]$.

	We will use the following notation.  Given $0 \le i_0 < \cdots
< i_p \le n$, let $U_{i_0\dots i_p} = U_{i_0}\cap\cdots \cap
U_{i_p}$.  Also, let $j_1 < \cdots < j_{n-p}$ be the complementary
indices, so that $\{i_0,\dots,i_p\}\cup\{j_1,\dots,j_{n-p}\} =
\{0,\dots,n\}$ is a disjoint union.  Finally, let
$\epsilon(i_0,\dots,i_p)$ be the sign of the permutation sending
$0,\dots,n$ to $i_0,\dots,i_p,j_1,\dots,j_{n-p}$ respectively.

	We now define some \v Cech cochains.  First, let $\gamma_p
\in C^p(\U,{\cal A}^{n,n-p})$ be given by
$$(\gamma_p)_{i_0\dots i_p} = \epsilon(i_0,\dots,i_p)\,
	\dbar\rho_{j_1}\wedge\cdots\wedge\dbar\rho_{j_{n-p}}
	\wedge\omega_g\ .$$
Note that each $\dbar\rho_j$ is divisible by $f_j$, so
that $\dbar\rho_{j_1}\wedge\cdots\wedge\dbar\rho_{j_{n-p}}$ is divisible
by $f_{j_1}\cdots f_{j_{n-p}}$.  This guarantees that $\gamma_p \in
C^p(\U,{\cal A}^{n,n-p})$.  Second, let $\widetilde\gamma_p \in
C^p(\U,{\cal A}^{n,n-p-1})$ be given by
$$(\widetilde\gamma_p)_{i_0\dots i_p} = \epsilon(i_0,\dots,i_p)\,
	\textstyle{\sum_{s=1}^{n-p}}
	(-1)^{s-1}\rho_{j_s}\,\dbar\rho_{j_1}\wedge\cdots\wedge
	\widehat{\dbar\rho_{j_s}}\wedge\cdots\wedge\dbar\rho_{j_{n-p}}
	\wedge\omega_g\ .$$
As above, we have $\widetilde\gamma_p \in C^p(\U,{\cal A}^{n,n-p-1})$.

	It is easy to see that
$$\dbar\widetilde\gamma_p = (n-p)\,\gamma_p\ ,\leqno(9)$$
and, if $\delta$ is the coboundary in the \v Cech complex, then
$$\delta \widetilde\gamma_p = (-1)^{p+1}\gamma_{p+1}\leqno(10)$$
(we will omit the straightforward but cumbersome proof).

	From (9) and (10), we get a class $[\gamma_p] \in H^p(X,{\cal
Z}^{n,n-p})$, and then a standard diagram chase, also using (9) and
(10), tells us that
$$\delta_p((n-p)[\gamma_p]) = (-1)^{p+1}[\gamma_{p+1}] \in
	H^{p+1}(X,{\cal Z}^{n,n-p-1})\ .$$
It follows that $n!\,[\gamma_0]$ maps to $(-1)^{n(n+1)/2}[\gamma_n]$
under the Dolbeault isomorphism .  However, $\gamma_n =
(\gamma_n)_{0\dots n} = \omega_g$, and using $\sum_{i=0}^n \rho_i =
1$, one can verify that $\gamma_0$ is given by
$$\eqalign{\gamma_0 &= \dbar\rho_1\wedge\cdots\wedge\dbar\rho_n\wedge
	\omega_g\cr &= \left({f_0\cdots f_n \textstyle{\sum_{i=0}^n}
	(-1)^i \overline{f}_i\, d\overline{f}_0\wedge\cdots\wedge
	\widehat{d\overline{f}_i}\wedge\cdots\wedge d\overline{f}_n
	\over (|f_0|^2 + \cdots + |f_n|^2)^{n+1}}\right)\wedge
	\omega_g \ ,\cr}$$
where the second equality follows from [GH, p.~655].  From here,
the proposition follows immediately.\qed
\medskip

	By Proposition A.1 in the appendix, we know that ${\rm
Tr}_X([\omega_g]) = \big({-1\over2\pi i}\big)^n\int_X\eta_g$.  Hence we
obtain the following integral representations of the toric residue.

\proclaim Theorem 6.3. If $X$ is simplicial and $f_0,\dots,f_n \in
S_\beta$ satisfy (3), then for $g \in S_\rho$, we have
$$\res(\omega_g) = {(-1)^{n(n-1)/2}n!\over(2\pi i)^n} \int_X
	{g\,\textstyle{\sum_{i=0}^n} (-1)^i
	\overline{f}_i\, d\overline{f}_0\wedge\cdots\wedge
	\widehat{d\overline{f}_i}\wedge\cdots\wedge d\overline{f}_n
	\wedge\Omega \over (|f_0|^2 + \cdots + |f_n|^2)^{n+1}}\ .$$
Furthermore, if $J$ is the toric Jacobian of $f_0,\dots,f_n$, then
$$\res(\omega_g) = {(-1)^{n(n-1)/2}n!\over(2\pi i)^n} \int_X
	{g\,\overline{J}\,\overline{\Omega}\wedge\Omega \over (|f_0|^2
	+ \cdots + |f_n|^2)^{n+1}}\ .\eqno\square$$

	Before giving our second set of integral formulas, we review
some symplectic geometry.  As above, $X$ is a projective simiplicial
toric variety, so that $X$ is a geometric quotient $(\C^{\Sigma(1)} -
Z)/G$.  Now let $G_\R = {\rm Hom}_\Z(A_{n-1}(X),S^1)$ be the compact
Lie group associated to $G$.  This group has Lie algebra ${\rm
Lie}(G_\R) = {\rm Hom}_\Z(A_{n-1}(X),\R)$.  The action of $G_\R$ on
$\C^{\Sigma(1)}$ is Hamiltonian, which gives the {\it moment map}
$$ \mu : \C^{\Sigma(1)} \longrightarrow {\rm Lie}(G_\R)^*$$
(see [A, Chapter VI, \S3.1]).  For us, the key property of $\mu$ is
that if $\xi \in {\rm Lie}(G_\R)^*$ is a regular value of $\mu$, then
$G_\R$ acts on $\mu^{-1}(\xi)$, and there is a natural isomorphism
$$ \mu^{-1}(\xi)/G_\R \simeq X$$
(see [A, Chapter VI, Proposition 3.1.1]).

	When $X = \P^n$, the moment map is $\mu(x_0,\dots,x_n) =
{1\over2}(|x_0|^2+\dots+|x_n|^2)$, so that $\mu^{-1}(\xi)$ is the
sphere $S^{2n+1}_r$ of radius $r = \sqrt{2\xi}$, and the map
$\mu^{-1}(\xi) \to \P^n$ is the Hopf fibration.  Another example is $X
= \P^1\times\P^1$, where $\mu(x,y;z,w) =
{1\over2}(|x|^2+|y|^2,|z|^2+|w|^2)$ and $\mu^{-1}(\xi_1,\xi_2) =
S^3_{r_1}\times S^3_{r_2}$, where $r_i = \sqrt{2\xi_i}$.  In general,
$\mu^{-1}(\xi)$ is a real manifold of dimension $|\Sigma(1)| + n$, and
we call the map
$$ \mu^{-1}(\xi) \longrightarrow X$$
the {\it generalized Hopf fibration\/} of $X$.  When $X$ is smooth,
this map is a genuine fibration with fiber $G_\R$, but in the
simplicial case, this is only true generically.

	Using the generalized Hopf fibration, we get two more integral
formulas for the toric residue.

\proclaim Theorem 6.4. If $X$ is simplicial and $f_0,\dots,f_n \in
S_\beta$ satisfy (3), let $d{\bf x} = \wedge_\rho dx_\rho$ for some
ordering of the $x_\rho$'s.  Then $\mu^{-1}(\xi)$ can be oriented so
that for $g \in S_\rho$, we have
$$\res(\omega_g) = {n!\over(2\pi i)^{|\Sigma(1)|}}
	\int_{\mu^{-1}(\xi)} {g\,\textstyle{\sum_{i=0}^n} (-1)^i
	\overline{f}_i\, d\overline{f}_0\wedge\cdots\wedge
	\widehat{d\overline{f}_i}\wedge\cdots\wedge d\overline{f}_n
	\wedge d{\bf x} \over (|f_0|^2 + \cdots + |f_n|^2)^{n+1}}\ .$$
Furthermore, if $J$ is the toric Jacobian of $f_0,\dots,f_n$, then
$$\res(\omega_g) = {n!\over(2\pi i)^{|\Sigma(1)|}}  \int_{\mu^{-1}(\xi)}
	{g\,\overline{J}\,\overline{\Omega}\wedge d{\bf x} \over (|f_0|^2
	+ \cdots + |f_n|^2)^{n+1}}\ .$$

\prf The first step is to relate ${\rm Lie}(G_\R)$ to the Euler
vector fields considered earlier.  Each $\vartheta \in {\rm
Lie}(G_\R)$ comes from a relation $\sum_\rho b_\rho n_\rho = 0$ where
$b_\rho \in \R$.  The Euler vector field $\theta = \sum_\rho b_\rho
x_\rho \partial /\partial x_\rho$ is a holomorphic vector field
tangent to $G$-orbits, and we also have a real vector field
$\vartheta$ tangent to $G_\R$ orbits.  These are connected by the
formula
$$ \vartheta = i(\theta - \overline{\theta})\ .$$
To see why, let $t$ be a real parameter.  If $x_\rho = u_\rho + i
v_\rho$ gives real and imaginary parts, then at the point $(x_\rho)
\in \C^{\Sigma(1)}$, we have
$$\vartheta = {d\over dt} (e^{ib_\rho t}(u_\rho+i v_\rho))_{|t=0} =
	(b_\rho(-v_\rho + i u_\rho))\ .$$
Using real coordinates, $\vartheta = \sum_\rho b_\rho\big(-v_\rho
\partial/\partial u_\rho + u_\rho \partial/\partial v_\rho\big)$,
which implies $\vartheta = i(\theta-\overline{\theta})$.

	Now suppose that $\vartheta_1,\dots,\vartheta_r$ form a basis of
${\rm Lie}(G_\R)$, where $r = |\Sigma(1)| - n$.  The corresponding
Euler vector fields are $\theta_1,\dots,\theta_r$, and we claim that
$$(\vartheta_1\wedge\cdots\wedge\vartheta_r)\ip(\overline{\Omega}\wedge
	d{\bf x}) = c\,i^r(-1)^{nr}\,\overline{\Omega}\wedge\Omega
	\ ,\leqno(11)$$
where $c$ satsifies $(\theta_1\wedge\cdots\wedge\theta_r)\ip d{\bf x}
= c\,\Omega$ as in Lemma 6.2.  To prove this, first note
that $\theta_j\ip\Omega = 0$ by Lemma 6.2, and thus
$\vartheta_j\ip\overline{\Omega} =
i(\theta_j-\overline{\theta}_j)\ip\overline{\Omega} = 0$.  Hence,
$$\vartheta_j\ip(\overline{\Omega}\wedge d{\bf x}) = (-1)^n
	\overline{\Omega} \wedge (\vartheta_j \ip d{\bf x}) = i (-1)^n
	\overline{\Omega}\wedge(\theta_j\ip d{\bf x})\ ,$$
and (11) now follows easily from Lemma 6.2.

	We next express ``integration over the fiber'' in terms of
interior multiplication.

\proclaim Lemma 6.5. Let $\eta$ be a $G_\R$-invariant form on
$\mu^{-1}(\xi)$ of degree $|\Sigma(1)| + n$, and let
$\vartheta_1,\dots,\vartheta_r$, $r = |\Sigma(1)|-n$, be a basis of
${\rm Lie}(G_\R)$.  Then we have the product formula
$$\int_{\mu^{-1}(\xi)} \eta = \int_X (\vartheta_1\wedge\cdots\wedge
	\vartheta_r)\ip \eta\cdot \int_{G_\R} \vartheta_1^*\wedge\cdots
	\wedge\vartheta_r^*\ ,$$
where $\vartheta_1^*,\dots,\vartheta_r^*$ is the dual basis of
$\vartheta_1,\dots,\vartheta_r$, regarded as invariant $1$-forms on
$G_\R$.

\prf Since we can remove sets of measure zero, we can replace
$\mu^{-1}(\xi)$ with a $G_\R$-stable open set $W$ where $G_\R$ acts
freely with quotient $X_0 \subset X$.  Then, by a partition of unity
argument, we can replace $X_0$ with an open set $U$ such that the
fibration is diffeomorphic to a product over $U$.  We can also assume
that $U$ has local coordinates $u_1,\dots,u_{2n}$.  Since our product
formula is invariant under diffeomorphism, we can replace the total
space by $U\times G_\R$.  For a small open set $V \subset G_\R$, we
can find local coordinates $t_1,\dots,t_r$ such that $\vartheta_i =
\partial /\partial t_i$.  Then
$$\eta = f(u_1,\dots,u_n)\,dt_1\wedge\cdots \wedge dt_r\wedge
	du_1\wedge\cdots \wedge du_{2n}$$
since $\eta$ is invariant under $G_\R$.  Hence the product formula
holds on $U\times V$, and another partition of unity
argument shows that it also holds on $U\times G_\R$.\qed
\medskip

	If we apply Lemma 6.5 to the $G$-invariant form
$$ \eta = {g\,\overline{J}\,\overline{\Omega}\wedge d{\bf x} \over
	(|f_0|^2 + \cdots + |f_n|^2)^{n+1}}$$
and use (11), we obtain
$$\int_{\mu^{-1}(\xi)} \eta = c\,i^r(-1)^{nr}\int_X {g\,\overline{J}\,
	\overline{\Omega}\wedge\Omega \over (|f_0|^2 + \cdots +
	|f_n|^2)^{n+1}} \cdot \int_{G_\R} \vartheta_1^*\wedge\cdots
	\wedge\vartheta_r^*\ .$$
If we can prove that
$$c\,\int_{G_\R} \vartheta_1^*\wedge\cdots\wedge\vartheta_r^* =
	\pm(2\pi)^r \ ,\leqno(12)$$
then the theorem will follow immediately from Theorem 6.3, provided we
adjust the orientation of $\mu^{-1}(\xi)$ to make the sign disappear.

	We will prove (12) using coordinates to compute the integral
explicitly.  First pick a subset $I = \{\rho_1,\dots,\rho_n\} \subset
\Sigma(1)$ such that the $n_{\rho_i}$ are linearly independent.  If
$\widehat{I}$ is the complement of $I$ in $\Sigma(1)$, then we have an
exact sequence
$$ 0 \longrightarrow \Z\otimes\widehat{I} \longrightarrow A_{n-1}(X)
	\longrightarrow F \longrightarrow 0\ ,\leqno(13)$$
where $F$ is a finite group.  Applying ${\rm Hom}_\Z(-,\R)$, we get
$\R^r = {\rm Hom}_\Z(\Z\otimes\widehat{I},\R) \simeq {\rm Lie}(G_\R)$,
and composing with the exponential map ${\rm Lie}(G_\R) \to
G_\R^\circ$, we get a covering map $\R^r \to G_\R^\circ$.  Here,
$G_\R^\circ \subset G_\R$ is the connected component of the identity.
For later purposes, note
$$[G_\R:G_\R^\circ] = |A_{n-1}(X)_{\rm tor}|\ .\leqno(14)$$

	Each vector field $\vartheta_i$ on $G_\R$ is determined by a
relation $\sum_\rho b_\rho^i n_\rho = 0$.  If $t_\rho$, $\rho \notin
I$, are the obvious coordinates on $\R^r$, the pull-back of
$\vartheta_i$ to $\R^r$ is $\sum_{\rho \notin I} b_\rho
\partial/\partial t_\rho$.  Thus $\vartheta_1\wedge\cdots
\wedge\vartheta_r$ pulls back to $\det(\hat{b}_I)\wedge_{\rho \notin
I} \partial/\partial t_\rho$, where $\det(\hat{b}_I)$ is as in
Lemma~6.2, and it follows that
$\vartheta_1^*\wedge\cdots\wedge\vartheta_r^*$ pulls back to
$\det(\hat{b}_I)^{-1} d{\bf t}$, where $d{\bf t} = \wedge_{\rho \notin
I} dt_\rho$.

	We have an exact sequence
$$ 0 \longrightarrow {\rm Hom}_\Z(A_{n-1}(X),2\pi\Z) \longrightarrow
	{\rm Lie}(G_\R) \longrightarrow G_\R^\circ \longrightarrow 0
	\ ,$$
and under ${\rm Lie}(G_\R) \simeq \R^r$, the lattice ${\rm
Hom}_\Z(A_{n-1}(X),2\pi\Z)$ maps to a sublattice $L' \subset
(2\pi\Z)^r = {\rm Hom}_\Z(\Z\otimes\widehat{I},2\pi\Z)$.  Applying
${\rm Hom}_\Z(-,\Z)$ to (13) gives the exact sequence
$$ 0 \longrightarrow L' \longrightarrow (2\pi\Z)^r \longrightarrow {\rm
	Ext}^1_\Z(F,\Z) \longrightarrow {\rm Ext}^1_\Z(A_{n-1}(X),\Z)
	\longrightarrow 0\ ,$$
which implies
$$ [(2\pi\Z)^r:L'] = {|F| \over |A_{n-1}(X)_{\rm tor}|}
	\ .\leqno(15)$$
Finally, the map $\Z\otimes\widehat{I} \to A_{n-1}(X)$ of (13) embeds in
a commutative diagram
$$\matrix{& & & & \Z\otimes\widehat{I} & = &
	\Z\otimes \widehat{I} & & \cr
	&& && \downarrow && \downarrow &&\cr
	0 & \to & M & \to & \Z\otimes\Sigma(1) & \to &
	A_{n-1}(X) & \to & 0\cr}$$
where the bottom row is exact (see [F, \S3.4]).  Since $\Z\otimes I =
{\rm coker}(\Z\otimes\widehat{I} \to \Z\otimes\Sigma(1))$, the snake
lemma and (13) imply
$$ |F| = |{\rm coker}(M \to \Z\otimes I)| = |\det(n_I)|\ ,\leqno(16)$$
where $\det(n_I)$ is as in the definition of $\Omega$.

	Using (14), (15) and (16), we see that $\int_{G_\R}
\vartheta_1^*\wedge\cdots\wedge\vartheta_r^*$ equals
$$\eqalign{|A_{n-1}(X)_{\rm tor}| \int_{G_\R^\circ}
	\vartheta_1^*\wedge\cdots\wedge\vartheta_r^* &=
	\pm|A_{n-1}(X)_{\rm tor}| \int_{R^r/L'} \det(\hat{b}_I)^{-1}
	d{\bf t}\cr &=
	\pm{|A_{n-1}(X)_{\rm tor}| \over\det(\hat{b}_I)\,[(2\pi\Z)^r:L']}
	\int_{\R^r/(2\pi\Z)^r} d{\bf t}\cr &=
	\pm{(2\pi)^r|F| \over\det(\hat{b}_I)} = \pm{(2\pi)^r|\det(n_I)|
	\over\det(\hat{b}_I)} \ , \cr}$$
where the sign $\pm1$ depends on whether $\R^r \to G_\R^\circ$ is
orientation preserving or not.  Hence
$$c\int_{G_\R} \vartheta_1^*\wedge\cdots\wedge\vartheta_r^* =
	\pm c\,{(2\pi)^r|\det(n_I)|
	\over\det(\hat{b}_I)} = \pm(2\pi)^r\,{c\,\det(n_I)
	\over\det(\hat{b}_I)} = \pm(2\pi)^r\ ,$$
where the last equality holds by (8).  This completes the proof of the
theorem. \qed
\medskip

	The simplest example is where $X= \P^n$.  Here, we've seen
that $\mu^{-1}(\xi)$ is a sphere $S^{2n+1}_r$, and then Theorem 6.4
tells us that $\res(\omega_g)$ is given by the integral
$${n!\over(2\pi i)^{n+1}}
	\int_{S_r^{2n+1}} {g\,\textstyle{\sum_{i=0}^n} (-1)^i
	\overline{f}_i\, d\overline{f}_0\wedge\cdots\wedge
	\widehat{d\overline{f}_i}\wedge\cdots\wedge d\overline{f}_n
	\wedge dx_0\wedge\cdots\wedge dx_n
	\over (|f_0|^2 + \cdots + |f_n|^2)^{n+1}}\ .$$
This formula appears in [GH, Chapter 5] and [T, \S5].

	We should remark that the theory of toric residues, as given
above, is not quite complete: one still needs to investigate whether
the Grothendieck local residue, as defined in (1), generalizes to the
toric case.  It would be very interesting to have such a formula.
\medskip

	{\bf Appendix. The Dolbeault Isomorphism and the Trace Map.}
In this appendix, let $X$ be a variety which is a compact orbifold (or
$V$-manifold).  As in \S6, this gives a Dolbeault isomorphism
$H_{\dbar}^{n,n}(X) \simeq H^n(X,\Omega^n_X)$.  Such a variety is also
Cohen--Macaulay, so that we have maps
$$\eqalign{{\rm Tr}_X : H^n(X,\Omega^n_X) &\to \C\ ,\quad\hbox{the
	trace map}\cr
	\textstyle{\int_X} : H_{\dbar}^{n,n}(X) &\to \C\
	,\quad\hbox{integration of $(n,n)$-forms\ .}\cr}$$
These maps are related as follows.

\proclaim Proposition A.1. Let $X$ be a compact orbifold variety, and
suppose that the $(n,n)$-form $\eta$ corresponds to $[\omega] \in
H^n(X,\Omega^n_X)$ under the Dolbeault isomorphism.  Then
$$ {\rm Tr}_X([\omega]) = \left(-1\over2\pi i\right)^n \int_X \eta\ .$$

\prf Consider the map $T : H^n(X,\Omega_X^n) \to \C$ defined by
$T([\omega]) = \int_X \eta$, where $\eta$ corresponds to $[\omega]$
under the Dolbeault isomorphism.  This map has all
of the formal properties of the trace map, except for the
normalization condition, which says it takes the value $1$ on the class
$$ \big[\omega_1\big] = \left[{\Omega_{\P^n} \over x_0\cdots
	x_n}\right] \in H^n(\P^n,\Omega^n_{\P^n})\ .$$
It follows that $T$ and ${\rm Tr}_X$ agree up to a constant, and
furthermore, the normalization condition implies the constant is
$T([\omega_1])$.  By Proposition 6.1, $T([\omega_1]) =
\int_{\P^n}\eta_1$, where
$$\eta_1 = (-1)^{n(n+1)/2} n!\, {\overline{\Omega}_{\P^n}\wedge
	\Omega_{\P^n} \over (|x_0|^2 + \cdots + |x_n|^2)^{n+1}}\ .$$
Hence, to prove the proposition, it suffices to show that $\int_{\P^n}
\eta_1 = (-2\pi i)^n$.

	If $\C^n \subset \P^n$ is the open set where $x_0 \ne 0$, then
$\int_{\P^n} \eta_1 = \int_{\C^n} \eta_1$.  On this open set, we can use
$t_i = x_i/x_0$ as coordinates.  Since $\Omega_{\P^n} = x_0^{n+1}
dt_1\wedge\cdots\wedge dt_n$, we have
$$\int_{\P^n}\eta_1 = (-1)^{n(n+1)/2}n!\,\int_{\C^n}
	{d\overline{t}_1\wedge\cdots\wedge
	d\overline{t}_n\wedge dt_1\wedge\cdots\wedge dt_n \over
	(1 + |t_1|^2 + \cdots + |t_n|^2)^{n+1}}\ .$$
Denoting real and imaginary parts by $t_j = u_j + iv_j$, we have
$d\overline{t}_j\wedge dt_j = 2i\,du_j\wedge dv_j$.  Thus
$d\overline{t}_1\wedge\cdots\wedge d\overline{t}_n\wedge
dt_1\wedge\cdots\wedge dt_n = (2i)^n (-1)^{n(n-1)/2}\,du_1\wedge dv_1
\wedge \cdots\wedge du_n\wedge dv_n$, so that
$$\int_{\P^n}\eta_1 = (-2i)^n n!\,\int_{\C^n} {du_1\wedge dv_1 \wedge
	\cdots\wedge du_n\wedge dv_n \over (1 + u_1^2 + v_1^2 + \cdots
	+ u_n^2 + v_n^2)^{n+1}}\ .$$
Since $du_1\wedge dv_1 \wedge \cdots\wedge du_n\wedge dv_n$ is the
volume form on $\C^n$, the above integral is a standard multiple
integral.  Using polar coordinates for each $u_j$ and $v_j$, the
multiple integral equals $\pi^n/n!$, and it follows that
$\int_{\P^n}\eta_1 = (-2\pi i)^n$.  This completes the proof.\qed
\medskip

	For the reader curious about the signs in Propositions 6.1 and
A.1, we suggest checking ${\rm Tr}_{\P^1}([\omega_1]) = {-1\over2\pi
i}\int_{\P^1}\eta_1 = 1$ in detail.
\medskip

\noindent {\bf References}
\medskip

\item{[A]} M.~Audin, {\sl The Topology of Torus Actions on Symplectic
Manifolds\/}, Progress in Math.~{\bf 93}, Birkh\"auser, Boston Basel
Berlin, 1991.
\smallskip

\item{[Bai]} W.~Baily, {\it The decomposition theorem for
$V$-manifolds\/}, Amer.~J.~Math.~{\bf 78} (1956), 862--888.
\smallskip

\item{[Bat]} V.~Batyrev, {\it Variations of the mixed Hodge structure of
affine hypersurfaces in algebraic tori\/}, Duke Math.~J.~{\bf 69}
(1993), 349--409.
\smallskip

\item{[BC]} V.~Batyrev and D.~Cox, {\it On the Hodge structure of
projective hypersurfaces in toric varieties\/}, Duke Math.~J.~{\bf 75}
(1994), 293--338.
\smallskip

\item{[BH]} W.~Bruns and J.~Herzog, {\sl Cohen--Macaulay Rings\/},
Cambridge Univ. Press, Cambridge, 1993.
\smallskip

\item{[CCD]} E.~Cattani, D.~Cox and A.~Dickenstein, {\it Residues in
toric varieties\/}, in preparation.
\smallskip

\item{[C]} D.~Cox, {\it The homogeneous coordinate ring of a toric
variety\/}, J.~Algebraic Geom., to appear.
\smallskip

\item{[D]} V.~Danilov, {\it The geometry of toric varieties\/},
Russian Math.~Surveys {\bf 33} (1978), 97--154.
\smallskip

\item{[F]} W.~Fulton, {\sl Introduction to Toric Varieties\/},
Princeton Univ.~Press, Princeton, 1993.
\smallskip

\item{[G]} P.~Griffiths, {\it On the periods of certain rational
integrals, I\/}, Ann.~of Math.~{\bf 90} (1969), 460--495.
\smallskip

\item{[GH]} P.~Griffiths and J.~Harris, {\sl Principles of Algebraic
Geometry\/}, John Wiley \& Sons, New York, 1978
\smallskip

\item{[GKZ]} I.~Gelfand, M.~Kapranov and A.~Zelevinsky, {\sl
Discriminants, Resultants, and Multidimensional Determinants\/},
Birkh\"auser, Boston Basel Berlin, 1994.
\smallskip

\item{[GW]} S.~Goto and K.~Watanabe, {\it On graded rings, I\/},
J.~Math.~Soc.~Japan~{\bf 30} (1978), 179--213.
\smallskip

\item{[Ha]} R.~Hartshorne, {\sl Residues and Duality\/}, Lecture Notes
in Math.~{\bf 20}, Springer-Verlag, Berlin Heidelberg New York, 1966.
\medskip

\item{[Ho]} M.~Hochster, {\it Rings of invariants of tori,
Cohen--Macaulay rings generated by monomials, and polytopes\/},
Ann.~of Math.~{\bf 96} (1972), 318--337.
\smallskip

\item{[K]} S.~Kleiman, {\it Toward a numerical theory of ampleness\/},
Ann.~of Math.~{\bf 84} (1966), 293--344.
\smallskip

\item{[O]} T.~Oda, {\sl Convex Bodies and Algebraic Geometry\/},
Springer-Verlag, Berlin Heidelberg New York, 1988.
\smallskip

\item{[PS]} C.~Peters and J.~Steenbrink, {\it Infinitesimal variation
of Hodge structure and the generic Torelli theorem for projective
hypersurfaces\/}, in {\sl Classification of Algebraic and Analytic
Manifolds} (K.~Ueno, editor), Progress in Math.~{\bf 39},
Birkh\"auser, Boston Basel Berlin, 1983, 399--463.
\smallskip

\item{[S]} S.~Sternberg, {\sl Lectures on Differential Geometry\/},
Prentice--Hall, Englewood Cliffs, 1964.
\smallskip

\item{[T]} A.~Tsikh, {\sl Multidimensional Residues and Their
Applications\/}, AMS, Providence, 1992.

\end

\item{[GHV]} W.~Greub, S.~Halperin and R.~Vanstone, {\sl Connections,
Curvature and Cohomology, Volume I}, Academic Press, New York, 1972.
\smallskip

We first describe orientations on $\mu^{-1}(\xi)$ and $G_\R$.
{}From the exact sequence
$$0 \longrightarrow M \longrightarrow \Z\otimes\Sigma(1)
	\longrightarrow A_{n-1}(X)\longrightarrow 0\leqno(12)$$
(see [F, \S3.4]), we obtain the fibration $G_\R \to {\rm
Hom}_\Z(\Z\otimes\Sigma(1),S^1) \to {\rm Hom}_\Z(M,S^1)$.  Then $d{\bf
x}$ determines an orientation of the total space and $\Omega$
determines an orientation of the base (by means of the basis
$m_1,\dots,m_n$ of $M$ used in the construction of $\Omega$).  This
gives an orientation of $G_\R$ where we put the base variables
before the fiber variables (this is the convention used in [GHV,
p.~288]).  Then the Hopf fibration $\mu^{-1}(\xi) \to X$ gives an
orientation on $\mu^{-1}(\xi)$ using the usual orientation of the
complex manifold $X$ and again putting base variables first.

	To keep track of the various orientations, first observe that
if we apply ${\rm Hom}_\Z(-,\R)$ to (18), we get a map of fibrations
$$\matrix{\R^r &\to& {\rm Hom}_\Z(\Z\otimes\Sigma(1),\R) &\to& {\rm
	Hom}_\Z(\Z\otimes I,\R)\cr \uparrow&&\|&&\uparrow\cr
	{\rm Lie}(G_\R) &\to& {\rm Hom}_\Z(\Z\otimes\Sigma(1),\R)
	&\to& N_\R\ .\cr}$$
The bottom fibration is oriented the same way that $G_\R$ was, so that
the exponential map ${\rm Lie}(G_\R) \to G_\R^\circ$ is orientation
preserving.  For the top fibration, we orient the fiber using $d{\bf
t}$, using the order induced from the $\rho$'s, and we orient the
base using $\rho_1,\dots,\rho_n$.  Hence we orient the total space
putting the base variables corresponding to $\rho_1,\dots,\rho_n$ at
the end, but otherwise preserving the order of the $\rho$'s.  With
these orientations, the vertical map to the right has determinant
$\det(n_I)$, and the ``identity'' map in the middle has determinant
$(-1)^I$ (this is the notation of Lemma 4.1).

	It follows that the sign of $(-1)^I\det(n_I)$ tells us whether
or not $\R^r \to G_\R^\circ$ is orientation preserving.  Thus

\prf We can assume $k = 0$ since the sign $(-1)^k$ comes from the
orientation of the local residue.  In light of Theorem 5.1, it
suffices to prove the proposition when $g \in \langle
f_0,\dots,f_n\rangle$ and when $g$ is the toric Jacobian of
$f_0,\dots,f_n$.

	If $g = \sum_{i=0}^n A_i\,f_i$, then near $x$, we can write
$g/f_0 = A_0 + \sum_{i=1}^n (A_i/f_0)\,f_i$, and standard properties
of residues (page 650 of [GH]) imply that
$$\res_{0,x} = \res_x\left({A_0\,\Omega \over f_1\cdots f_n}\right)\
	.$$
However, the meromorphic $n$-form $A_0\Omega/(f_1\cdots f_n)$ has
polar divisor contained in $D_1\cup\cdots\cup D_n$, and it follows
from the Global Residue Theorem (page 656 of [GH]) that the sum over
all $x \in D_1\cap\cdots\cap D_n$ is zero.  Since $g \in \langle
f_0,\dots,f_n\rangle$ implies $\res(\omega_g) = 0$
by Proposition~2.3, we get equality in this case.

	Next, suppose that $g$ is the toric Jacobian $J$.  In this
case, we claim that each local residue is a local intersection
product, i.e., for each $x \in D_1\cap\cdots\cap D_n$, we have
$$ \res_{0,x}(\omega_J) = (D_1\cdots D_n)_x\ .$$
Once this is proved, the sum of the local residues is the intersection
product $D_1\cdots D_n = D_1^n$ (remember that $D_1,\dots,D_n$ are
linearly equivalent).  In the proof of Theorem 5.1, we saw that
$\res(\omega_J)$ also equals $D_1^n$.  Thus the theorem will follow
once we show that each $\res_{0,x}(\omega_J)$ is a local intersection
product.

	To prove this, we use the following argument of Cattani and
Dickenstein.  Given $x$ as above, we can find a $n$-dimensional cone
$\sigma$ such that $x$ lies in the affine toric variety $X_\sigma$.
We also know that $X_\sigma$ is an affine space since $X$ is smooth.
If the $1$-dimensional cones of $\sigma$ are $I =
\{\rho_1,\dots,\rho_n\}$, then the proof of Theorem 1.9 of [BC] shows
that $x_{\rho_1},\dots,x_{\rho_n}$ give coordinates for $X_\sigma$
once we set $x_\rho = 1$ for all $\rho \notin I$.

	If we let $\tilde f_i$ denote the polynomial obtained from
$f_i$ by setting $x_\rho = 1$ for $\rho \notin I$, then the toric
Jacobian becomes
$$ \tilde J = {\det\pmatrix{\tilde f_0 & \cdots & \tilde f_n\cr
	{\partial \tilde f_0/\partial x_{\rho_1}} & \cdots &
	{\partial \tilde f_n / \partial x_{\rho_1}}\cr
	\vdots & & \vdots \cr
	{\partial \tilde f_0 / \partial x_{\rho_n}} & \cdots &
	{\partial \tilde f_n / \partial x_{\rho_n}}\cr}
	\over \det(n_I)}\ .$$
Expanding along the top row, it follows that
$$ \tilde J \equiv {\tilde f_0 \det(\partial\tilde f_i/\partial
	x_{\rho_j}) \over \det(n_I)} \bmod \tilde f_1,\dots,\tilde
	f_n\ .$$
Note also that in these coordinates, $\Omega = \det(n_I) dx_{\rho_1}
\wedge \cdots \wedge dx_{\rho_n}$.  Combining these equations and
using properties of local residues, we see that
$$ \res_{0,x}(\omega_J) = \res_x\left({(\tilde J/\tilde f_0)\,\Omega
	\over \tilde f_1 \cdots \tilde f_n}\right) =
	\res_x\left({\det(\partial\tilde f_i/\partial
	x_{\rho_j}) \over \tilde f_1 \cdots \tilde f_n}\right)\ .$$
The latter residue is well-known to equal the local intersection
multiplicity $(D_1\cdots D_n)_x$ (page 663 of [GH]).  This completes
the proof of the proposition.\qed
\medskip

	When $X$ is a compact complex manifold and $D_0,\dots,D_n$ are
divisors with empty intersection, one could ask if the analog of
Proposition 5.2 is true for any meromorphic form $\omega$ on $X$ with
$D_0 \cup \cdots \cup D_n$ as polar set.  Cattani and Dickenstein have
proved that this is true under the addition hypothesis that the
divisors $D_0,\dots,\widehat{D_k},\dots,D_n$ are a global regular
sequence.